\documentclass[11pt, a4paper, twoside, leqno]{article}
\usepackage[T1]{fontenc}
\usepackage[utf8]{inputenc}
\usepackage{lmodern}
\usepackage{amsmath,amssymb,amsthm,mathtools}
\usepackage[margin=26mm]{geometry}
\usepackage{microtype}
\usepackage[hidelinks]{hyperref}
\usepackage{enumitem,booktabs,tabularx}
\usepackage{xcolor}
\usepackage{bookmark}
\numberwithin{equation}{section}
\newtheorem{theorem}{Theorem}[section]
\newtheorem{proposition}[theorem]{Proposition}
\newtheorem{lemma}[theorem]{Lemma}

\theoremstyle{remark}

\newtheorem{example}[theorem]{Example}
\newcommand{\R}{\mathbb R}
\newcommand{\N}{\mathbb N}
\newcommand{\Om}{\Omega}
\newcommand{\M}{\mathcal M}
\newcommand{\supp}{\operatorname{supp}}

\newcommand{\ord}{\operatorname{ord}}
\newcommand{\diag}{\operatorname{diag}}
\newcommand{\norm}[1]{\left\lVert#1\right\rVert}

\hypersetup{hidelinks}

\let\oldthanks\thanks
\renewcommand{\thanks}[1]{%
	\renewcommand{\thefootnote}{}
	\oldthanks{#1}
	\renewcommand{\thefootnote}{\arabic{footnote}}
}

\begin{document}

\title {\textbf{ $L^p(\mathbb{R}^3)$-estimates for maximal averages along mixed\\
		homogeneous hypersurfaces with vanishing Hessian }	
	\author{Wenjuan Li, \hspace{0.2cm} Huiju Wang}
	\thanks{2020 {\it Mathematical Subject Classification}: 42B25; 42B20.}
	\thanks{{\it Keywords and phrases:} Maximal operator, mixed homogeneous, vanishing Hessian}
	\thanks{This work is supported by the National Key R\&D Program of China (No.2023YFA1010800); Natural Science Foundation of China (No.12301113; No.12271435).}}

\date{}
\maketitle
\thispagestyle{empty}

\begin{abstract}
We study maximal averages over the surfaces parametrized by $c+\Phi$ in $\R^3$, where $c$ is a constant and $\Phi$ is mixed homogeneous and real analytic away from the origin. Under the assumption that the Hessian determinant of $\Phi$ vanishes identically, we classify the possible graphing functions for all positive weights. For surfaces in each resulting class, we establish sharp
$L^{p}$-estimates  for the corresponding maximal averages.
\end{abstract}

\tableofcontents

\section{Introduction}

A fundamental problem in the theory of maximal averages is to determine how the geometry of a surface affects the $L^p$-boundedness of the associated maximal operator. Besides curvature and
contact orders, the position of the
surface relative to the origin is also very important:
translating the surface changes the family
of averages, since the translation is
multiplied by the dilation parameter.
This distinction is particularly relevant
when an affine tangent plane passes through
the origin. The estimates available under
the usual transversality assumption need
not describe the behavior in that situation;
see \cite{IKM05,IKM10,Zimmermann14}.

This paper focuses primarily on the related problems in $\mathbb{R}^3$.
Mixed homogeneous surfaces provide a natural
setting in which to study the interaction
between these features. Ikromov--Kempe--M\"uller
\cite{IKM05} considered surface parameterized by
$c+\Phi$, where $\Phi$ is mixed homogeneous,
and established sharp estimates under
finite-order curvature hypotheses when
$c\ne0$. More recently,
Buschenhenke--Li--M\"uller--Wang \cite{BLMW}
studied the maximal operator associated with mixed homogeneous surfaces passing
through the origin under the assumption
that the Hessian determinant does not
vanish identically. Their results identify
the critical exponents in terms of the
degeneracies of that determinant.
The case
\[
 \det D^2\Phi\equiv0
\]
requires a separate description. It includes
the familiar polynomial cylinders, but,
when analyticity at the origin is not
assumed, it also includes non-polynomial
conical surface.

In this article we determine the admissible
phases and the sharp $L^p$ ranges
in this identically vanishing Hessian
class. We assume that $\Phi$ is real
analytic on the entire punctured plane
and mixed homogeneous of degree one,
with arbitrary positive weights.
In particular, the weights may coincide
or be greater than or equal to one.
Starting from the global classification
of G\'alvez--Nelli \cite{GN16}, we derive
the precise forms compatible with these
weights. Up to invertible linear ambient
transformations, the resulting surfaces
are planes, monomial cylinders, or
 cones. This also
identifies exactly when a non-polynomial
phase can occur.
Our main new result determines the sharp
$L^p$-estimates for maximal averages
along non-planar analytic cones which do not pass through the origin.

Let $\kappa=(\kappa_1,\kappa_2)$, where
$\kappa_1,\kappa_2>0$, and define
\begin{equation}\label{eq:dilations}
 \delta_r(x)
 =(r^{\kappa_1}x_1,r^{\kappa_2}x_2),
 \qquad
 |\kappa|=\kappa_1+\kappa_2,
 \qquad r>0.
\end{equation}
A function $\Phi$ on
$\Om:=\R^2\setminus\{0\}$ is called
$\kappa$-homogeneous of degree $\alpha$ if
\[
 \Phi(\delta_rx)=r^\alpha \Phi(x)
 \qquad(x\in\Om,\ r>0).
\]
We write
\[
 \Phi_i=\partial_i\Phi,\qquad
 \Phi_{ij}=\partial_i\partial_j\Phi,\qquad
 \mathrm{H}\Phi=\det D^2\Phi.
\]

When $\kappa_1=\kappa_2=1$, we say that
$\Phi$ is \emph{homogeneous of degree one}
in the usual sense, namely
\[
 \Phi(\lambda x)=\lambda\Phi(x)
 \qquad(x\in\Om,\ \lambda>0).
\]
We use the term \emph{positively homogeneous of degree one}
to distinguish this isotropic scaling
from mixed homogeneity.
Equivalently, in polar coordinates,
\[
 \Phi(\rho\omega_\theta)=\rho a(\theta),
 \qquad a(\theta)=\Phi(\omega_\theta).
\]

Throughout the article, we assume
\begin{equation}\label{eq:assumptions}
 \Phi\in C^\omega(\Om;\R),\qquad
 \Phi(\delta_rx)=r\Phi(x),\qquad
 \mathrm{H}\Phi\equiv0.
\end{equation}
No upper bound on the positive weights
is imposed, and the weights may coincide.
We do not assume analyticity at the origin.

Every $\Phi$ satisfying
\eqref{eq:assumptions} has a continuous
extension with $\Phi(0)=0$.
Indeed, for each $x\ne0$, there is a unique
$r=r(x)>0$ such that
\[
 r^{-2\kappa_1}x_1^2
 +r^{-2\kappa_2}x_2^2=1.
\]
Writing $x=\delta_r\omega$ with $|\omega|=1$,
we have
\[
 |\Phi(x)|=r|\Phi(\omega)|
 \le r\sup_{|\omega|=1}|\Phi(\omega)|.
\]
Since both weights are positive,
$r\to0$ as $x\to0$.

For an analytic function $\Phi$, the
order $\ord \Phi(x)$ is the smallest
nonnegative integer $j$ for which
some derivative of total order $j$
is nonzero at $x$. It is infinity
if all derivatives vanish there.
In particular, the order is zero
when $\Phi(x)\ne0$.

For $c\in\R$, consider the surface
\begin{equation}\label{eq:surface}
 S_c=
 \{(x_1,x_2,c+\Phi(x)):x\in\Om\}.
\end{equation}
Its closure contains $(0,0,c)$.
In particular, the closure contains
the ambient origin exactly when $c=0$.
When referring to a surface passing
through the origin, we include this
closure convention.

Fix a density satisfying
\begin{equation}\label{eq:density}
 \psi\in C_c^\infty(\R^2),\qquad
 \psi\ge0,\qquad \psi(0)>0.
\end{equation}
For $c\in \mathbb{R}$, we define the averaging operator along the surface $S_c$ and the corresponding maximal operator by
\begin{equation}\label{eq:averages}
 \begin{split}
 A_t^cf(y)
 &:=
 \int_{\R^2}
 f\bigl(y-t(x,c+\Phi(x))\bigr)\psi(x)\,dx,\\
 \M_cf(y)
 &:=\sup_{t>0}|A_t^cf(y)|.
 \end{split}
\end{equation}
Specially, when \(c=0\), i.e., when the surface passes through the origin, we denote the averaging and maximal operators by \(A_{t}^{0}\) and \(\mathcal{M}_{0}\), respectively.

\subsection{Classification under positive weights}\label{section1.1}

The following statement gives the
explicit mixed homogeneous forms
of the known cylindrical and conical
alternatives. Its geometric ingredient
is Proposition~\ref{quoted:GN};
the weight calculations are supplied
in Section~\ref{sec:classification}.

\begin{theorem}[Classification under positive weights]
\label{thm:classification}
Assume that $\Phi$ satisfies \eqref{eq:assumptions}, where
$\kappa_1,\kappa_2>0$.
Then the following alternatives hold.

\begin{enumerate}[label=\textup{(\roman*)},leftmargin=*]

\item Suppose $\kappa_1\ne\kappa_2$.
Then $\Phi$ has one of the forms
\[
 \Phi(x)=\ell(x)
 \quad\text{or}\quad
 \Phi(x)=\ell(x)+bx_j^m,
\]
where
\[
 \ell(x)=\ell_1x_1+\ell_2x_2,\qquad
 \ell_i=0\ \text{whenever }\kappa_i\ne1,
 \quad i=1,2.
\]
In the second case,
\begin{equation}\label{eq:unequalform}
 j\in\{1,2\},\qquad b\ne0,\qquad
 m\in\N,\quad m\ge2,\qquad m\kappa_j=1.
\end{equation}
The zero linear form is allowed.
Since the weights are distinct,
at most one coefficient of $\ell$
can be nonzero.

\item Suppose
$\kappa_1=\kappa_2\ne1$.
Then either $\Phi\equiv0$, or there
are an integer $m\ge2$, a constant
$b\ne0$, and a nonzero real linear
form $L$ such that
\begin{equation}\label{eq:equalform}
 \kappa_1=\kappa_2=\frac1m,\qquad
 \Phi(x)=bL(x)^m.
\end{equation}

\item Suppose
$\kappa_1=\kappa_2=1$.
Then the possibilities are exactly
\begin{equation}\label{eq:coneform}
 \Phi(\rho\omega_\theta)=\rho a(\theta),
 \qquad
 \omega_\theta=(\cos\theta,\sin\theta),
 \quad \rho>0,
\end{equation}
where $a$ is an arbitrary real analytic
$2\pi$-periodic function.
For every such function,
$\mathrm{H}\Phi\equiv0$ automatically.

\end{enumerate}

Conversely, every listed function
satisfies the stated analyticity,
homogeneity, and vanishing Hessian
determinant conditions with the
indicated weights.
\end{theorem}

The theorem has a simple consequence:
outside the weight pair $(1,1)$,
every admissible function is a
polynomial, although analyticity
at the origin was not assumed.
The non-polynomial possibilities
occur in the positively homogeneous of degree one
conical class.


For a non-planar cone, define
\begin{equation}\label{eq:angularmultiplicity}
 \begin{gathered}
 B_\Phi(\theta):=a(\theta)+a''(\theta),\\
 N_\Phi(\theta):=\ord B_\Phi(\theta),\qquad
 M_\Phi:=\max_{\theta\in[0,2\pi)}
          \bigl(N_\Phi(\theta)+2\bigr).
 \end{gathered}
\end{equation}
Here $B_\Phi$ is not identically zero:
otherwise $a''+a \equiv 0$ and $\Phi$ is linear.
Analyticity on the circle therefore
makes $M_\Phi$ a finite integer
at least two.

The factor $B_\Phi$ measures the
angular component of the Hessian.
Its order is different from the
order of $\mathrm{H}\Phi$, which vanishes
identically throughout this article.
Lemma~\ref{lem:contactcone} explains
why a zero of order $N_{\Phi}$ of $B_\Phi$
gives transverse tangent-plane
contact of order $N_{\Phi}+2$.

\subsection{Sharp $L^{p}$-estimates for maximal operators}

Based on the classifications stated in Theorem \ref{thm:classification}, we present  the following theorem, which gives the sharp $L^{p}$-estimates for maximal operators $\mathcal{M}_{c}$ defined by (\ref{eq:averages}).

\begin{theorem}
\label{thm:maximal}
For any $c\in \mathbb{R}$, $\M_c$ is bounded on $L^\infty(\R^3)$.
For $1\le p<\infty$, its boundedness
is characterized as follows.

\begin{enumerate}[label=\textup{(\roman*)},leftmargin=*]

\item If $\Phi$ is linear, then $\M_0$
is bounded on $L^p$ if and only
if $p>1$.
For $c\ne0$, the operator $\M_c$
is unbounded on every finite
$L^p$.

\item Suppose
\[
 \Phi(x)=\ell(x)+bL(x)^m,
\]
where $\ell$ is linear, $L$ is a
nonzero real linear form, $b\ne0$,
and $m\ge2$ is an integer.
Then $\M_0$ is bounded on $L^p$
if and only if $p>2$.
For $c\ne0$, the operator $\M_c$
is bounded on $L^p$ if and
only if $p>m$.

\item Suppose $\Phi$ is positively homogeneous of degree one and is
not linear.
Let $M_\Phi$ be defined by
\eqref{eq:angularmultiplicity}.
Then $\M_0$ is bounded on $L^p$
if and only if $p>1$.
For $c\ne0$, the operator $\M_c$
is bounded on $L^p$ if and
only if $p>M_\Phi$.

\end{enumerate}
\end{theorem}

 We write $p_{crit}$ for the
critical exponent.
For comparison, the conclusions are
summarized below in terms of the
critical exponent $p_{crit}$.
When $p_{crit}<\infty$, the sharp strong
range is $p_{crit}<p\le\infty$.
The notation $p_{crit}=\infty$ means
that boundedness holds only on
$L^\infty$.

\begin{center}
\begin{tabular}{@{}lcc@{}}
\toprule
Class of the surface
& $c=0$ & $c\ne0$\\
\midrule
Planar
& $p_{crit}=1$ & $p_{crit}=\infty$\\
Non-planar cylindrical, degree $m\ge2$
& $p_{crit}=2$ & $p_{crit}=m$\\
Non-planar conical
& $p_{crit}=1$ & $p_{crit}=M_\Phi$\\
\bottomrule
\end{tabular}
\end{center}


\subsection{Examples}

\begin{example}[A cone with planar rays]
\label{ex:angular}
Let
\[
 a(\theta)=1+\frac13\cos(2\theta),
 \qquad
 \Phi(\rho\omega_\theta)=\rho a(\theta).
\]
Then
\[
 B_\Phi(\theta)=a(\theta)+a''(\theta)
 =1-\cos(2\theta)=2\sin^2\theta.
\]
The zeros of $B_{\Phi}$ have order two, so
$M_\Phi=4$. According to Theorem \ref{thm:maximal}-(iii),
we have $p_{crit}=1$ when $c=0$, whereas
  $p_{crit}=4$ when $c\ne0$.

The rays $\theta=0,\pi$ consist
of planar points, i.e., where $D^{2}\Phi$ is the zero matrix, but the surface does not coincide with a plane in any neighborhood of those points.
\end{example}

\begin{example}\label{ex:weights}
For
\[
 (\kappa_1,\kappa_2)=(1,1/3),
 \qquad \Phi(x)=\ell_1x_1+bx_2^3,
 \qquad b\ne0,
\]
all assumptions in Theorem \ref{thm:maximal}-(ii) are satisfied.
The critical exponent is $p_{crit}=2$ when
$c=0$ and $p_{crit}=3$ when $c\ne0$.
\end{example}



\section{Auxiliary estimates}\label{sec:quoted}
In this section we collect several standard facts and known maximal estimates. They will be used in Section~\ref{sec:necessary} and in the boundedness arguments below. We include references for the reader's convenience. The classification argument in Section~\ref{sec:classification} is independent of the present section and can be read
separately.

\subsection{Linear changes of coordinates}\label{sec:reductions}

 For a parametrized patch $X: \Omega \rightarrow \mathbb{R}^{3}$ with
parameter density $\psi$, we  write
$\M_{X,\psi}$ for the
maximal operator defined by
\begin{equation}\label{MX}
\mathcal{M}_{X,\psi}f(y): =\sup_{t>0}\biggl|\int_{\mathbb{R}^{2}}f\bigl(y-t X(x)\bigl)\psi(x)dx\biggl|, \quad \quad y \in \mathbb{R}^{3}.
\end{equation}
If $A\in GL(3,\R)$ is an invertible matrix, let
$\widetilde f=f\circ A^{-1}$, then
\[
 \M_{X,\psi}f(A^{-1}y)
 =\M_{AX,\psi}\widetilde f(y).
\]
 For $1\le p<\infty$, \cite[Lemma~2.2.1]{Zimmermann14}  gives
\[\|\mathcal{M}_{X,\psi}\|_{L^{p} \rightarrow L^{p} } =\|\mathcal{M}_{AX,\psi}\|_{L^{p} \rightarrow L^{p} }.\]
The same conclusion for $p=\infty$
is immediate. In particular, for $(x,x_{3}) \in \mathbb{R}^{3}$, the transformation
\[
 (x,x_3)\longmapsto(x,x_3-\ell(x))
\]
removes the linear term in
Theorem~\ref{thm:classification}
without changing the constant  $c$ in the formula for $\mathcal{M}_{c}$ in (\ref{eq:averages}).
An invertible horizontal map with
first coordinate $L(x)$ reduces
$bL(x)^m$ to a coordinate monomial.

For an invertible matrix $B\in GL(2,\mathbb{R})$, under a linear change of variables in the parameter plane
$w=Bx$, the density $\psi$ becomes
\[
 \widetilde\psi(w)
 =\psi(B^{-1}w)|\det B|^{-1},
\]
which preserves the conditions
in \eqref{eq:density}.

Next, we recall some maximal estimates
which will be used in the proof of Theorem  \ref{thm:maximal}.

\subsection{A uniform small-parameter surface estimate}

Let $\mathcal U\subset\R^2$ be an open neighborhood of the point
$x^0\in\mathbb{R}^2$, and let
\[
 \varphi_{\rm{pr}}:\mathcal U\longrightarrow\R
\]
be smooth. Suppose that, for an
integer $k\ge2$,
\begin{equation}\label{pr}
 \partial_2^k\varphi_{\rm{pr}}(x^0)\ne0.
\end{equation}
For a sufficiently small smooth real-valued perturbation $\varphi_{r}$, put
$\varphi=\varphi_{\rm{pr}}+\varphi_{r}$.

For  for $y\in\R^3$, define
\[
 \begin{split}
 \mathcal N_{\varepsilon}^{\varphi}f(y)
 :=
 \sup_{t>0}\left|
 \int
 f\bigl(
 y-t(x,1+\varepsilon\varphi(x))
 \bigr)\psi(x)\,dx
 \right|,
 \end{split}
\]
where   $\varepsilon>0$ is
independent of the dilation
parameter $t$, and $\psi\in C_{c}^{\infty}(U)$.

\begin{proposition}
\label{quoted:epsilon} (\cite[Proposition~4.5]{IKM10}) Assume that $\phi_{pr}$ satisfy (\ref{pr}) and the neighborhood $\mathcal U$ of $x^0$ is sufficiently small. Then
for each $p>k$, there exist  an integer
$N$, positive numbers $\delta$ and
$\varepsilon_0$, and a constant
$C_p>0$ such that for $\varepsilon>0$ sufficiently small, the maximal operator $\mathcal N_{\varepsilon}^{\varphi}$ satisfies the following estimate
\[
 \norm{\mathcal N_{\varepsilon}^{\varphi}f}
       _{L^p}
 \le
 C_p\varepsilon^{-1/p}
 \norm f_{L^p}
\]
whenever
\[
 \| \varphi_{r} \|_{C^N(\mathcal U)}<\delta.
\]
The constant $C_{p}$ is independent
of $\varepsilon$ and of the admissible
perturbation $\varphi_{r}$.
\end{proposition}


\section{Classification and local geometry when
\texorpdfstring{$\mathrm{H}\Phi\equiv0$}
{the Hessian determinant vanishes}}
\label{sec:classification}

We first recall the global classification due to Gálvez–Nelli and derive the restrictions imposed by mixed homogeneity. We then prove Theorem~\ref{thm:classification}, treating the cases of distinct and equal weights separately. Finally, we describe the tangent planes and planar points of the classified surfaces; for cones, we relate the angular vanishing order of
$a +a^{\prime\prime}$ to the transverse contact order used in the maximal estimates.

\subsection{A global classification}

A surface is called cylindrical here if, after
a rotation in the horizontal variables, its
graphing function has the form
\[
 G(se+tv)=F(s)+\beta t,
\]
where $e$ and $v$ are fixed orthonormal vectors.
Indeed, varying $t$ traces straight lines
on the surface, all of which are  parallel to the fixed
vector $(v,\beta)\in\R^3$.

A surface is called conical with vertex at
$(0,0,C)$ if
\[
G(\rho\omega_\theta)=C+\rho a(\theta).
\]
For fixed $\theta$, varying $\rho>0$
traces a straight ray from that vertex.

\begin{proposition}(\cite[Theorem~1]{GN16})
\label{quoted:GN}
Let $G\in C^\omega(\R^2\setminus\{0\};\R)$
satisfy
\[
 G_{11}G_{22}-G_{12}^2=0.
\]
Then at least one of the following holds:
\begin{enumerate}[label=\textup{(\roman*)},leftmargin=*]
\item There are orthonormal vectors
$e,v\in\R^2$, a constant $\beta\in\R$,
and a function $F\in C^2(\R)$ such that
\[
G(se+tv)=F(s)+\beta t
 \qquad\text{whenever }se+tv\ne0.
\]
\item There are a constant $C\in\R$ and
a real analytic $2\pi$-periodic function
$a$ such that
\[
G(\rho\omega_\theta)=C+\rho a(\theta),
 \qquad \rho>0.
\]
\end{enumerate}
\end{proposition}

In alternative \textup{(i)}, the function
$F$ is in fact real analytic on $\R$.
To see this, evaluate the representation
on the line $se+v$, which never passes
through the origin:
\[
 F(s)=G(se+v)-\beta.
\]
The right-hand side is real analytic for
every $s\in\R$.

The alternatives need not be mutually exclusive:
a hyperplane can be represented both
cylindrically and conically.

\subsection{Consequences of mixed homogeneity}

Differentiating the homogeneity identity $\Phi(\delta_{r}x)=r\Phi(x)$
in the base variables and in the dilation
parameter gives
\begin{equation}\label{eq:derivativescaling}
 \Phi_i(\delta_rx)
 =r^{1-\kappa_i}\Phi_i(x),
 \qquad
 \Phi(x)
 =\kappa_1x_1\Phi_1(x)
  +\kappa_2x_2\Phi_2(x).
\end{equation}
We shall also need the Hessian scaling
identity. Put
\[
 D_r=\diag(r^{\kappa_1},r^{\kappa_2}),
 \qquad \delta_rx=D_rx.
\]
Differentiating twice gives
\[
 D_r^{\mathsf T}D^2\Phi(D_rx)D_r
 =rD^2\Phi(x).
\]
Since $D_r$ is diagonal,
\begin{equation}\label{eq:hessianscaling}
 D^2\Phi(D_rx)
 =rD_r^{-1}D^2\Phi(x)D_r^{-1}.
\end{equation}
In particular, Hessian rank is preserved
under mixed dilations.

Positive mixed homogeneity also gives
the continuous extension $\Phi(0)=0$,
as explained after \eqref{eq:assumptions}.
Therefore, if the conical alternative
in Proposition~\ref{quoted:GN}-(ii) holds,
the constant $C$ therein must be zero. Then we have
\[
 \Phi(\rho\omega_\theta)=\rho a(\theta).
\]
Such a function is positively homogeneous of degree one.
We next determine when a non-planar
  cone can also be $\kappa$-homogeneous of degree  one.

\begin{lemma}\label{lem:coneweights}
Suppose that $\Phi$ is real analytic
on $\Om$, positively homogeneous of degree one, and $\kappa$-homogeneous of degree
one for positive weights $\kappa_1,\kappa_2$.
If $\Phi$ is not linear, then
\[
 \kappa_1=\kappa_2=1.
\]
\end{lemma}

\begin{proof}
Suppose first that
$\kappa_1=\kappa_2$.
The positively homogeneity and mixed Euler identities
give
\[
 \Phi=x\cdot\nabla\Phi,\qquad
 \Phi=\kappa_1 x\cdot\nabla\Phi=\kappa_2 x\cdot\nabla\Phi.
\]
Thus $(1-\kappa_i)\Phi=0$, $i=1,2$.
Since a non-linear $\Phi$ is not the
zero function, $\kappa_1=\kappa_2=1$.

Now suppose $\kappa_1\ne\kappa_2$.
Solving the two Euler identities for
$x_1\Phi_1$ and $x_2\Phi_2$ gives
\begin{equation}\label{eq:twoeuler}
 \begin{gathered}
 x_1\Phi_1=A\Phi,\qquad
 x_2\Phi_2=B\Phi,\\
 A=\frac{1-\kappa_2}{\kappa_1-\kappa_2},
 \qquad
 B=\frac{\kappa_1-1}{\kappa_1-\kappa_2},
 \qquad A+B=1.
 \end{gathered}
\end{equation}

By Taylor  expansion of $\Phi$ near the point $(0,1)$, we obtain
\[
 \Phi(x_1,x_2)
 =\sum_{n=0}^\infty a_n(x_2)x_1^n.
\]
The first equation in
\eqref{eq:twoeuler} implies
\[
 (n-A)a_n(x_2)=0
 \qquad(n\ge0).
\]
At least one coefficient $a_{n}(x_{2})$ is nonzero:
otherwise $\Phi$ would vanish on an
open set and hence on all of $\Om$.
Therefore $A$ is a non-negative integer.

Similarly, expansion in $x_2$ near
$(1,0)$ shows that $B$ is a non-negative
integer. Since $A+B=1$,
\[
 (A,B)=(1,0)
 \quad\text{or}\quad
 (A,B)=(0,1).
\]

If $(A,B)=(1,0)$, then $x_2\Phi_2=0$,
so $\Phi_2=0$ on $\Om$.
The remaining identity is
$\Phi=x_1\Phi_1$.
Differentiating it in $x_1$ gives
$x_1\Phi_{11}=0$, and hence
$\Phi_{11}=0$ everywhere.
Thus the Hessian matrix vanishes
identically and $\Phi$ is affine.
Positively degree-one homogeneity removes
the constant term, so $\Phi$ is linear.
The other case is identical with the
variables interchanged.
This contradicts the assumption that
$\Phi$ is not linear.
\end{proof}

Lemma \ref{lem:coneweights}   fails  if the analyticity of $\Phi$ is  assumed
only on a quadrant. For example, let $\Phi(x)=x^{3}_{1}x^{-2}_{2}$. Then $\Phi$
 is not only  positively homogeneous but also $\kappa$-homogeneous of degree one whenever
 $3\kappa_{1} - 2\kappa_{2} =1$.

\subsection{Classification for distinct weights}

\begin{proof}[Proof of Theorem~\ref{thm:classification},
part \textup{(i)}]
Assume $\kappa_1\ne\kappa_2$.

We first dispose of the planar case.
If $D^2\Phi\equiv0$, then connectedness
of $\Om$ gives
\[
 \Phi(x)=\ell_1x_1+\ell_2x_2+d.
\]
The continuous extension $\Phi(0)=0$
implies $d=0$.
The derivative scaling identity gives
\[
 \ell_i=r^{1-\kappa_i}\ell_i
 \qquad(r>0).
\]
Consequently,
\[
 \ell_i=0\quad\text{unless }\kappa_i=1.
\]
This is the stated linear alternative.

Suppose from now on that $\Phi$ is not linear.
By the standing assumption
\[
\mathrm{H}\Phi
=
\Phi_{11}\Phi_{22}
-
\Phi_{12}^{\,2}
\equiv0
\quad\text{on }\Omega,
\]
Proposition~\ref{quoted:GN}
therefore applies and gives either a
cylindrical or a conical representation
of $\Phi$.
The conical alternative is excluded by
Lemma~\ref{lem:coneweights}, since
$\kappa_1\ne\kappa_2$.
Hence
\[
\Phi(x)=F(e\cdot x)+\beta\,v\cdot x,
\]
where $e,v$ are orthonormal and
$F$ is real analytic on $\mathbb R$.

Differentiating this representation
gives
\[
 D^2\Phi(x)
 =F''(e\cdot x)\,ee^{\mathsf T}.
\]
Because $\Phi$ is not linear, $F''$
is not identically zero.
Choose $x$ with $F''(e\cdot x)\ne0$.
At every point where the Hessian is
nonzero, its kernel is the same line:
\[
 \ker D^2\Phi=\R v.
\]
We note that here we used the fact that
$v$ is orthogonal to $e$.

By \eqref{eq:hessianscaling}, the Hessian
matrix of $\Phi$ at $D_rx$ is also nonzero, and
\[
 \ker D^2\Phi(D_rx)
 =D_r\ker D^2\Phi(x).
\]
Both sides must equal the fixed kernel
line $\R v$. Hence $D_rv$ is parallel
to $v$ for every $r>0$.
Taking $r=2$ gives
\[
 0=\det(v,D_2v)
 =v_1v_2
   \bigl(2^{\kappa_2}-2^{\kappa_1}\bigr).
\]
Since the weights are distinct,
$v_1v_2=0$.
Thus the fixed flat direction $v$ is
parallel to a coordinate axis.

Suppose first that $v$ is parallel
to the vector $(0,1)$. Absorbing possible signs
into the one-variable function and
the constant, we may write
\[
 \Phi(x)=F_0(x_1)+\beta_0x_2,
 \qquad F_0\in C^\omega(\R), \quad \quad \beta_{0} \in \mathbb{R}.
\]
Since $\Phi_2=\beta_0$, the homogeneity of $\Phi_{2}$ implies
\[
 \beta_0=0\quad\text{unless }\kappa_2=1.
\]
Thus $\beta_0x_2$ is   $\kappa$-homogeneous of degree one. Since $\Phi$ is $\kappa$-homogeneous of degree one,
$F_{0}$ is also $\kappa$-homogeneous
of degree one. Euler's identity yields
\[
 F_0(x_{1})=\kappa_1 x_{1}F_0'(x_{1}).
\]
The Taylor expansion of $F_{0}(x_{1})$ at $x_{1}=0$ gives
\[
 F_0(x_{1})=\sum_{n=0}^\infty a_nx_{1}^n.
\]
Comparison of coefficients implies
\[
 (1-\kappa_1n)a_n=0
 \qquad(n\ge0).
\]
The function $F_0$ is not identically
zero, since otherwise $\Phi$ would
be linear. There is therefore a
positive integer $m$ such that
\[
 m\kappa_1=1,
\]
and the only nonzero coefficient is
$a_m=b$, where $b$ is some nonzero real constant.
Thus $F_0(x_{1})=bx_{1}^m$ near zero and then everywhere,
since $F_{0}$ is real analytic   on
$\R$.
The case $m=1$ would make $\Phi$
linear, so $m\ge2$. We obtain
\[
 \Phi(x)=\beta_0x_2+bx_1^m.
\]

If the flat direction $v$ is parallel
to the vector $(1,0)$, the same argument gives
\[
 \Phi(x)=\beta_0x_1+bx_2^m,
 \qquad m\kappa_2=1.
\]
Together with the planar case,
these are exactly the alternatives
in part \textup{(i)}.

\end{proof}

\subsection{Classification for equal weights}

\begin{proof}[Proof of Theorem~\ref{thm:classification},
parts \textup{(ii)} and \textup{(iii)}]
Suppose
$\kappa_1=\kappa_2=\sigma>0$
and set $d=1/\sigma$.
The homogeneity of $\Phi$ gives
\begin{equation}\label{eq:ordinary}
 \begin{gathered}
 \Phi(\lambda x)=\lambda^d\Phi(x),
 \qquad \lambda>0.
 \end{gathered}
\end{equation}

First assume $\sigma\ne1$, or
equivalently, $d\ne1$.
A linear function satisfying (\ref{eq:ordinary}) must be identically zero.
If $\Phi$ is nonzero, it therefore  must be
non-linear.

Proposition~\ref{quoted:GN}
gives either a cylindrical  or a conical
representation   for the function $\Phi$.
The conical representation would
imply
\[
\Phi(\lambda x)=\lambda\Phi(x),
\qquad \lambda>0.
\]
Comparing this with
$\Phi(\lambda x)=\lambda^d\Phi(x)$,
where $d\ne1$ and $\Phi\not\equiv0$,
excludes the conical case.
Hence
\[
\Phi(se+tv)=F(s)+\beta t,
\]
where $e,v$ are orthonormal and $F$ is
real analytic on $\mathbb R$.
By (\ref{eq:ordinary}), we have
\[
 F(\lambda s)+\beta\lambda t
 =\lambda^dF(s)+\beta\lambda^dt.
\]
For fixed $\lambda$ and $s$, comparison
of the coefficients of $t$ gives
\[
 \beta(\lambda-\lambda^d)=0.
\]
Since $d\ne1$, we can choose
$\lambda>0$ with $\lambda\ne\lambda^d$.
Thus $\beta=0$.

The Euler identity now reduces to
\[
 sF'(s)=dF(s).
\]
Then, by arguments similar to those in the proof of part (i), we obtain
\[
 F(s)=bs^m,\qquad b\ne0, \quad s\in \mathbb{R}.
\]
Thus
\[
 \sigma=\frac1m,\qquad
 \Phi(x)=b(e\cdot x)^m.
\]
Replacing $e$ by a general nonzero linear
form $L$ only changes the constant
$b$, hence this proves part (ii).

 Now assume $\sigma=1$.
 According to   (\ref{eq:ordinary}),
$\Phi$ is determined by its restriction to
the unit circle. Parametrizing the circle
by $\omega_\theta=(\cos\theta,\sin\theta)$, we
write
\[
a(\theta)=\Phi(\omega_\theta).
\]
Then
\[
\Phi(\rho\omega_\theta)=\rho a(\theta),
\qquad \rho>0.
\]
The function $a$ is real analytic and
$2\pi$-periodic. Conversely, the Hessian
condition imposes no further restriction
on $a$: Euler's identity for a
 homogeneous function of degree one gives
$D^2\Phi(x)x=0$ for $x\ne0$, and hence
$\mathrm{H}\Phi\equiv0$. This proves part \textup{(iii)}.

\end{proof}

\subsection{Tangent planes and their remainders}

For the surface
\[
 X_c(x)=(x,c+\Phi(x)),
\]
the tangent plane at $X_c(x^0)$ has equation
\[
 z_3=c+\Phi(x^0)
     +\nabla\Phi(x^0)\cdot(z'-x^0),
 \qquad z'=(z_1,z_2).
\]
Its constant term is
\begin{equation}\label{eq:transversality}
 T_c(x^0)
 :=
 c+\Phi(x^0)-x^0\cdot\nabla\Phi(x^0).
\end{equation}
Thus the tangent plane contains the
ambient origin if and only if
$T_c(x^0)=0$.
We call the surface \emph{transversal}
at $X_c(x^0)$ when $T_c(x^0)\ne0$.
Adding a linear function to $\Phi$
does not change this quantity.

Define the tangent-plane remainder by
\begin{equation}\label{remaider}
 R_{x^0}(h)
 :=
 \Phi(x^0+h)-\Phi(x^0)
 -\nabla\Phi(x^0)\cdot h.
\end{equation}
It measures the vertical difference
between the surface and its tangent plane.


For an positively homogeneous cone,
Euler's identity gives
\begin{equation}\label{constant}
 x\cdot\nabla\Phi(x)=\Phi(x),
 \qquad\text{and hence}\qquad
 T_c(x)=c.
\end{equation}
Consequently, all tangent planes
contain the origin when $c=0$,
whereas none contains the origin
when $c\ne0$.
The next subsection determines the
order of the remainder in the
direction transverse to a radial
generator.

For completeness, a point $x^{0}$ over the surface is
called \emph{planar} when both principal
curvatures vanish there.
For a surface, this is equivalent to
$D^2\Phi(x^0)=0$.
It is a pointwise condition and does
not imply that the surface agrees with
a plane on a neighborhood.

\subsection{Angular curvature and contact for cones}

Let
\[
\Phi(\rho\omega_\theta)=\rho a(\theta),
\qquad
\omega_\theta=(\cos\theta,\sin\theta),
\qquad
e_\theta=(-\sin\theta,\cos\theta).
\]
Here $e_\theta$ is the unit angular
direction, and
\[
\omega_\theta'=e_\theta,
\qquad
e_\theta'=-\omega_\theta.
\]

Differentiating
$\Phi(\rho\omega_\theta)=\rho a(\theta)$
with respect to $\rho$ and $\theta$ gives
\begin{equation}\label{gradientformula}
\nabla\Phi(\rho\omega_\theta)
\cdot\omega_\theta=a(\theta),
\qquad
\nabla\Phi(\rho\omega_\theta)
\cdot e_\theta=a'(\theta).
\end{equation}
Since $\omega_\theta,e_\theta$ form an
orthonormal basis of $\R^2$, we obtain
\[
\nabla\Phi(\rho\omega_\theta)
=
a(\theta)\omega_\theta
+
a'(\theta)e_\theta.
\]
This gradient is constant along each
radial ray. Its derivative in the
$e_\theta$ direction is
\[
D^2\Phi(\rho\omega_\theta)e_\theta
=
\frac{1}{\rho}
\frac{d}{d\theta}
\bigl(a(\theta)\omega_\theta
+a'(\theta)e_\theta\bigr)
=
\frac{a(\theta)+a''(\theta)}{\rho}
e_\theta.
\]
Consequently, writing
\[
B_\Phi(\theta):=a(\theta)+a''(\theta),
\]
we have
\begin{equation}\label{eq:conehessian}
	D^2\Phi(\rho\omega_\theta)
	=
	\frac{B_\Phi(\theta)}{\rho}
	e_\theta e_\theta^{\mathsf T},
	\qquad
	T_c(\rho\omega_\theta)=c.
\end{equation}
The last identity follows from (\ref{constant}).

If $B_\Phi\equiv0$, then
$a''+a=0$, so
\[
 a(\theta)=A\cos\theta+B\sin\theta
\]
and $\Phi(x)=Ax_1+Bx_2$.
Thus, for a non-planar cone,
$B_\Phi$ is not identically zero.

A nonzero real analytic function
on the circle has finitely many
zeros, each of finite order.
With the convention that the order
is zero at a nonzero value, recall
\[
 N_\Phi(\theta)=\ord B_\Phi(\theta),
 \qquad
 M_\Phi=\max_\theta
          \bigl(N_\Phi(\theta)+2\bigr).
\]
In particular, $M_\Phi$ is finite
and at least two. The following computation identifies
this number with a transverse
tangent-plane contact order.

\begin{lemma}\label{lem:contactcone}
Let $\Phi$ define a non-planar cone,
and fix
\[
 x^0=\rho_0\omega_{\theta_0},
 \qquad \rho_0>0.
\]
Use horizontal coordinates
\begin{equation}\label{formularofx}
 x=x^0+u\omega_{\theta_0}+v e_{\theta_0}
   =(\rho_0+u)\omega_{\theta_0}
      +v e_{\theta_0}.
\end{equation}
If
\[
 k=N_\Phi(\theta_0)+2,
\]
then
\[
 R_{x^0}(u\omega_{\theta_0}+v e_{\theta_0})
 =v^kU(u,v),
\]
where $U$ is real analytic near
$(0,0)$ and $U(0,0)\ne0$, and $R_{x_0}$ is defined by (\ref{remaider}).
\end{lemma}

\begin{proof}

Restricting to $|u|<\rho_0/2$, we have
$\rho_0+u>0$. Set
\[
q=\frac{v}{\rho_0+u}.
\]
Then  we   rewrite the formula  (\ref{formularofx}) as
\[
x=(\rho_0+u)
\bigl(\omega_{\theta_0}
+qe_{\theta_0}\bigr).
\]
Euler's identity gives
$\Phi(x^0)=\nabla\Phi(x^{0})\cdot x^0$, and the homogeneity of $\Phi$ gives
$\nabla\Phi(x^{0}) = \nabla\Phi(\omega_{\theta_{0}})$.
Thus,  we have
\[
 \begin{split}
 R_{x^0}(u\omega_{\theta_0}+v e_{\theta_0})
 &=\Phi(x)-\nabla\Phi(x^{0})\cdot x = \Phi(x)-\nabla\Phi(\omega_{\theta_{0}})\cdot x \\
 & = (\rho_{0} + u) \bigl[\Phi\bigl(\omega_{\theta_0}
+qe_{\theta_0}\bigr)- \nabla\Phi(\omega_{\theta_{0}}) \cdot \bigl(\omega_{\theta_0}
+qe_{\theta_0}\bigr)   \bigl].
 \end{split}
\]

In order to describe the expression $\bigl[\Phi\bigl(\omega_{\theta_0}
+qe_{\theta_0}\bigr)- \nabla\Phi(\omega_{\theta_{0}}) \cdot \bigl(\omega_{\theta_0}
+qe_{\theta_0}\bigr)   \bigl]$, we define
\[
 b_0(\theta):=a(\theta)-\nabla\Phi(\omega_{\theta_{0}})\cdot\omega_\theta.
\]
The gradient formula in
\eqref{gradientformula} gives
\[
 b_0(\theta_0)=0,\qquad
 b_0'(\theta_0)=0.
\]
Also,
\[
 b_0''+b_0=B_\Phi.
\]
Set $N=N_\Phi(\theta_0)$.
Then
\[
 B_\Phi^{(j)}(\theta_0)=0
 \quad(0\le j<N),\qquad
 B_\Phi^{(N)}(\theta_0)\ne0.
\]
Successively differentiating the
identity for $b_0$ yields
\[
 b_0^{(j+2)}(\theta_0)
 =
 B_\Phi^{(j)}(\theta_0)
 -b_0^{(j)}(\theta_0).
\]
Starting from $b_0=b_0'=0$ at
$\theta_0$, we obtain
\[
 b_0^{(j)}(\theta_0)=0
 \quad(0\le j<N+2),
\]
and
\[
 b_0^{(N+2)}(\theta_0)
 =B_\Phi^{(N)}(\theta_0)\ne0.
\]
Thus $b_0$ has a zero of exact
order $k=N+2$.

Since $b_0$ has a zero of exact
order $k$ at $\theta_0$, there is a
real analytic function $b_*$ such that
\[
b_0(\theta_0+\delta)
=
\delta^k b_*(\delta),
\qquad
b_*(0)
=
\frac{b_0^{(k)}(\theta_0)}{k!}
\ne0.
\]
Also, for $q$ near $0$,
\[
\omega_{\theta_0}+qe_{\theta_0}
=
\sqrt{1+q^2}\,
\omega_{\theta_0+\arctan q}.
\]
By homogeneity,
\[
\begin{aligned}
	\Phi(\omega_{\theta_0}+qe_{\theta_0})
	-\nabla\Phi(\omega_{0})\cdot(\omega_{\theta_0}+qe_{\theta_0}) &=\sqrt{1+q^2} [\Phi(\omega_{\theta_0+\arctan q})
	-\nabla\Phi(\omega_{0})\cdot\omega_{\theta_0+\arctan q} ] \\
	&=
	\sqrt{1+q^2}\,
	b_0(\theta_0+\arctan q)
	=
	q^k b_1(q),
\end{aligned}
\]
where
\[
b_1(q)
=
\sqrt{1+q^2}
\left(\frac{\arctan q}{q}\right)^k
b_*(\arctan q).
\]
The quotient $(\arctan q)/q$ extends
analytically to $q=0$ with value one.
Hence $b_1$ is real analytic and
$b_1(0)=b_*(0)\ne0$.

Consequently,
\[
 \begin{split}
 R_{x^0}(u\omega_{\theta_0}+v e_{\theta_0})
  =(\rho_0+u)q^kb_1(q)=v^k(\rho_0+u)^{1-k}
   b_1\left(\frac{v}{\rho_0+u}\right).
 \end{split}
\]
Therefore
\[
 U(u,v)
 =(\rho_0+u)^{1-k}
   b_1\left(\frac{v}{\rho_0+u}\right)
\]
is analytic near $(0,0)$ and satisfies
\[
 U(0,0)=\rho_0^{1-k}b_1(0)\ne0.
\]
\end{proof}

Thus the remainder vanishes identically
in the radial direction $v=0$.
In the transverse direction, its first
nonzero term has degree
$N_\Phi(\theta_0)+2$.
At a non-planar point this degree is
two; at a planar point it is greater
than two but remains finite.
The largest transverse contact order is
exactly $M_\Phi$.


\section{Necessary conditions}\label{sec:necessary}
We recall the notation
\[
 X_c(x)=(x,c+\Phi(x)),\qquad
 A_t^cf(y)=\int_{\R^2}f(y-tX_c(x))\psi(x)\,dx,
 \qquad
 \M_cf(y)=\sup_{t>0}|A_t^cf(y)|.
\]
By the arguments in Subsection~\ref{sec:reductions}, we can reduce the
functions classified in Theorem~\ref{thm:classification}
to planes, monomial cylinders, and cones, without changing the $L^{p}$-norm of the
corresponding maximal operator except by a constant factor. In this section,
we briefly explain how the necessary conditions for the $L^{p}$-boundedness of
maximal operators associated with planes, monomial cylinders, and cones can be derived.

\subsection{The hyperplane integrability condition}



The following  theorem is a consequence of \cite[Theorem 2]{IS96}.

\begin{lemma}\label{lem:contacttest}
	Suppose that   $ \psi  \ge J_{0}> 0$  on the region $I_0\times(-w_0,w_0)$, where $I_0$ is an interval of positive length and $w_0>0$.
	Suppose further that there exist
	$\eta\in\mathbb R^3\setminus\{0\}$,
	$d\in\mathbb R\setminus\{0\}$, an integer $k\ge2$, and a
	continuous real-valued function
	\[
	U:I_0\times(-w_0,w_0)\longrightarrow\mathbb R
	\]
	such that
	\begin{equation}\label{eq:contacttest}
		\eta\cdot X_{c}(x)-d=x_{2}^kU(x),
		\qquad
		0<C_1\le |U(x)|\le C_2
	\end{equation}
	uniformly for each $x =(x_{1},x_{2}) \in I_0\times(-w_0,w_0)$. Then $\mathcal{M}_{c}$ is unbounded on $L^p$ whenever
	$1<p\le k$.
\end{lemma}

\begin{proof}
	Consider the affine plane
	\[
	H=\{z\in\mathbb R^3:\eta\cdot z=d\}.
	\]
	It does not contain the origin, since $d\ne0$. Because $\eta$
	has not been normalized, the distance from $z$ to $H$ is
	\[
	\operatorname{dist}(z,H)
	=\frac{|\eta\cdot z-d|}{|\eta|}.
	\]
	In particular, \eqref{eq:contacttest} implies that
	\[
	\operatorname{dist}(X_{c}(x),H)
	=\frac{|x_{2}|^k|U(x)|}{|\eta|}
	\le \frac{C_2}{|\eta|}|x_{2}|^k.
	\]
	
	The   assumption on the function $\psi$ therefore yields, for $1<p\le k$,
	\[
	\begin{aligned}
		\int_{\mathbb R^2}
		\operatorname{dist}(X_{c}(x),H)^{-1/p} \psi(x)\,dx
		&\ge
		J_0\int_{I_0}\int_{-w_0}^{w_0}
		\operatorname{dist}(X_{c}(x),H)^{-1/p}
		\,dx\\
		&\ge
		J_0\left(\frac{|\eta|}{C_2}\right)^{1/p}
		|I_0|\int_{-w_0}^{w_0}|x_{2}|^{-k/p}\,dx_{2}
		=\infty.
	\end{aligned}
	\]
	Here the last integral diverges because $k/p\ge1$. Since $H$ does not
	contain the origin, the necessary condition in
	 \cite[Theorem 2]{IS96} rules out $L^p$-boundedness for $1<p\le k$.
\end{proof}

\subsection{Translated planes and cylinders}\label{necessaryPC}

Now we consider the case where $S_{c}$ is a translated plane.
By a change of variables, we may remove the linear term of $\Phi$,
so that $S_c$ takes the form
\[
 X_c(x)=(x_1,x_2,c).
\]
If $c\ne0$, take $H=\{z_3=c\}$.
The distance from the surface to $H$ vanishes
on every patch of positive measure, so
\cite[Theorem 2]{IS96}
rules out $L^p$-boundedness of $\mathcal{M}_{c}$ for every finite $p>1$.

For a cylindrical surface, by non-degenerate linear transformations, we may just consider the normal form
\[
X_c(x)=(x_{1},x_{2},c+bx_{1}^m),
\qquad b\ne0,\quad m\ge2.
\]
Assume $c\ne0$ and take $\eta=(0,0,1)$. Then
\[
\eta\cdot X_c(x)-c=bx_{1}^m.
\]
Since $\psi(0)>0$, the parameter density is bounded away from zero
on a sufficiently small rectangle about $(0,0)$.
Lemma~4.2 applies with   $k=m$. It follows that
the maximal operator $\mathcal{M}_{c}$ is unbounded on $L^{p}$ for $1<p\le m$.

\subsection{Translated cones}\label{subsection4.3}

Let $\Phi$ define a non-planar cone. Recall that
\[
B_\Phi=a+a'',\qquad
N_\Phi(\theta)=\operatorname{ord}_\theta B_\Phi(\theta),
\qquad
M_\Phi=\max_\theta\bigl(N_\Phi(\theta)+2\bigr).
\]
Assume $c\ne0$, and choose $\theta_0$ such that
\[
N_\Phi(\theta_0)+2=M_\Phi.
\]
Since $\psi$ is positive near the origin, we may choose
$\rho_0>0$ sufficiently small so that
$x^0=\rho_0\omega_{\theta_0}$ satisfies $\psi(x^{0})>0$.

Using the horizontal coordinates
\[
x(u,v)=x^0+u\omega_{\theta_0}+v e_{\theta_0},
\]
 we have
\[
X_{c}(u,v)
=
\bigl(x(u,v),\,c+\Phi(x(u,v))\bigr).
\]
 Set
\[
\eta=(-\nabla\Phi(x^{0}),1).
\]
Euler's identity and
Lemma \ref{lem:contactcone} give
\[
\eta\cdot X_{c}(u,v)-c
=
\Phi(x(u,v))-\nabla\Phi(x^{0})\cdot x(u,v)
=R_{x_0}(u\omega_{\theta_0}+v e_{\theta_0})
=
v^{M_\Phi}U(u,v),
\qquad U(0,0)\ne0,
\]
where $U$ is real analytic near $(0,0)$, and $R_{x_0}$ is defined by (\ref{remaider}).

By continuity, we can restrict $(u,v)$ to a rectangle
$I_0\times(-v_0,v_0)$ on which
$0<C_1\le |U(u,v)|\le C_2$. The horizontal coordinate
change is orthogonal, and $\psi(x^0)>0$; after shrinking
the rectangle if necessary, the parameter density is
bounded away from zero by a positive constant there. Applying
Lemma~4.2 with   $k=M_\Phi$
shows that the maximal operator $\mathcal{M}_{c}$ is unbounded on $L^{p}$ for
$1<p\le M_\Phi$.

\subsection{A standard failure of
\texorpdfstring{$L^1$}{L1}-boundedness}

\begin{lemma}\label{lem:Lone}
The maximal operator $\mathcal{M}_{c}$
fails to be bounded on $L^1 $.
\end{lemma}

\begin{proof}
This is a standard fact in harmonic analysis, see \cite[Chapter I, Section 3]{Stein}.
\end{proof}

\section{Surfaces passing through the origin}
\label{sec:origin}

In this section, we show the $L^{p}$-estimates for maximal operators associated with planes, cylinders, and cones that pass through the origin. Throughout this section, we mainly study the $L^p$-estimates for $\mathcal{M}_{0}$.

\subsection{Planes and cylinders}

We first consider the case when $\Phi$ is linear.
After removing the linear term,
we may assume that $\Phi=0$.
If $\supp\psi\subset[-R,R]^2$, then
\[
 \begin{split}
 |A_t^0f(y)|
 &\le\norm\psi_\infty\,t^{-2}
 \int_{[-Rt,Rt]^2}
 |f(y_1-w_1,y_2-w_2,y_3)|\,dw\\
 &\le4R^2\norm\psi_\infty\,
 H^{(2)}[f(\,\cdot\,,\,\cdot\,,y_3)](y_1,y_2),
 \end{split}
\]
where $ H^{(2)}$ is the two-dimensional Hardy–Littlewood maximal operator. The
boundedness of $ H^{(2)}$ on $L^{p}(\mathbb{R}^{2})$, together with  Fubini's
theorem, gives the $L^{p}(\mathbb{R}^{3})$-boundedness of $\mathcal{M}_{0}$ for $p>1$.
Lemma~\ref{lem:Lone} yields that  $\mathcal{M}_{0}$ is unbounded on
$L^1(\mathbb{R}^{3})$.

Next, we consider the cylindrical case.
For
\[
 \Phi(x)=\ell(x)+bL(x)^m, \quad \quad b\neq 0, m \ge 2,
\]
by   non-degenerate linear transformations, it suffices to consider the case where
\[
 \Phi(x)=bx_1^m,\qquad b\ne0,\quad m\ge2.
\]
Then the sufficiency of the condition  $p>2$ for the  $L^{p}$-boundedness of $\mathcal{M}_{0}$
follows from
\cite[Theorem~1.8]{Li16} by the first author. The necessity of the  condition $p>2$ follows from Lemma 6.3.1 in \cite{Zimmermann14}.

\subsection{Cones}

For
\[
 \Phi(\rho\omega_\theta)=\rho a(\theta),
\]
put
\[
 V(\theta)=(\omega_\theta,a(\theta)).
\]
Note that $V(\theta) \neq 0$ for each $\theta$.
For a nonzero vector $V\in\R^3$,
define
\[
 H_V^+f(y)
 =
 \sup_{h>0}\frac1h
 \int_0^h|f(y-sV)|\,ds.
\]
By rotational invariance, we may fix the direction of $V$,
while its length is absorbed into the parameter $h$.
The one-dimensional Hardy--Littlewood
maximal theorem gives
\begin{equation}\label{Vmaximal}
 \norm{H_V^+f}_{L^{p}}\le C_p\norm f_{L^{p}},
 \qquad p>1,
\end{equation}
uniformly in $V\ne0$.

If $\supp\psi\subset\{|x|\le R\}$,
polar coordinates yield
\[
 A_t^0f(y)
 =
 \int_0^{2\pi}\int_0^R
 f(y-t\rho V(\theta))
 \psi(\rho\omega_\theta)\rho\,d\rho\,d\theta.
\]
For fixed $t$ and $\theta$,
\[
 \begin{split}
 \int_0^R|f(y-t\rho V(\theta))|\rho\,d\rho
 &\le\frac Rt
 \int_0^{tR}|f(y-sV(\theta))|\,ds\\
 &\le R^2H_{V(\theta)}^+f(y).
 \end{split}
\]
Hence
\begin{equation}\label{eq:coneradial}
 \M_0f(y)
 \le R^2\norm\psi_\infty
 \int_0^{2\pi}H_{V(\theta)}^+f(y)\,d\theta.
\end{equation}
Minkowski's inequality and (\ref{Vmaximal}) yield
\[
 \norm{\M_0f}_{L^{p}}
 \le C_p\norm f_{L^{p}},
 \qquad p>1.
\]
Lemma~\ref{lem:Lone} implies that $p >1$ is also necessary for $\mathcal{M}_{0}$ to be $L^{p}$-bounded.

\section{Translated surfaces}
\label{sec:translated}

\subsection{Translated planes and cylinders}

For a translated plane, the necessary
condition in Subsection~\ref{necessaryPC}
excludes every finite exponent.
The $L^\infty$ bound is immediate.

For the translated cylindrical case, the sufficient part follows directly from
weighted estimate for translated
finite-type hypersurfaces in
\cite[Theorem 1.22]{LWZ23}.
The necessary condition in
Subsection~\ref{necessaryPC}
excludes $p\le m$.

\subsection{Translated cones}
\label{subsec:translatedcones}


Finally, we consider the $L^{p}$-boundedness of $\mathcal{M}_{c}$ associated with translated cones.
The necessary condition $p>M_{\Phi}$ follows from
Subsection~\ref{subsection4.3}.
To prove the sufficient condition, it suffices to establish the $L^{p}$-boundedness of $\mathcal{M}_{c}$ for any fixed exponent $p$ with $M_\Phi<p<\infty$.

We first prove the following proposition.

\begin{proposition}[Uniform estimate on a fixed annulus]
\label{prop:cone-fixed-annulus}
Let $\Phi$ be as in (\ref{eq:coneform}).
We choose a non-negative function
$\tilde{\psi}\in C_c^\infty(\mathbb R^2\setminus\{0\})$ supported in
the annulus $\{x\in \mathbb{R}^{2}: c_{1}<|x| < c_{2}\}$, where $0<c_{1}<c_{2}$ are fixed constants.
For $\varepsilon > 0$, we define
\[
 \widetilde{\mathcal M}  f(y)
: =
 \sup_{t>0}\left|
 \int_{\mathbb R^2}
 f\bigl(y-t(x,c+\varepsilon\Phi(x))\bigr)
 \tilde{\psi}(x)\,dx
 \right|.
\]
Then for each $p>M_{\Phi}$, we have the following $L^{p}$-estimates for  $\widetilde{\mathcal{M}}$.

\textbf{(1)} When $\varepsilon =1$, we have
\begin{equation}\label{away0}
\|\widetilde{\mathcal{M}} \|_{L^{p}  \rightarrow L^{p} } \le C_{p}.
\end{equation}

\textbf{(2) } For $\varepsilon$ sufficiently small, we have
\begin{equation}\label{eq:cone-fixed-annulus}
 \|\widetilde{\mathcal M} \|_{L^p \rightarrow L^p }
 \le
 C_{p} \varepsilon^{-\frac{1}{p}}.
\end{equation}
\end{proposition}

\begin{proof}
\textbf{(1)}  We will apply  \cite[Theorem~1.2]{IKM10} to prove part (1). For this purpose,
we need to show that the transversality and finite-type assumptions of
\cite[Theorem~1.2]{IKM10} are satisfied.

 Fix $x^0=\rho_0\omega_{\theta_0}\in\operatorname{supp}\psi$,
and let
\[
 k=2+\operatorname{ord}_{\theta_0}B_\Phi
 \le M_\Phi.
\]
We first compute the height at
\[
 P^0
   =(x^0,c+ \Phi(x^0)).
\]
We use the orthogonal
coordinates
\[
 x=x^0+u\omega_{\theta_0}+v e_{\theta_0}.
\]
By Lemma \ref{lem:contactcone},
\begin{equation}\label{eq:cone-local-remainder}
 \Phi(x)-\Phi(x^{0})-\nabla\Phi(x^{0})\cdot (x-x^{0}) = R_{x^0}(u\omega_{\theta_0}+v e_{\theta_0})=v^kU(u,v),
 \qquad U(0,0)\ne0,
\end{equation}
where $U$ is real analytic near $(0,0)$. Euler's identity
gives $\Phi(x^0)=\nabla\Phi(x^{0})\cdot x^0$, so the affine tangent plane
to the surface $S_{c}=\{(x,c+\Phi(x)): x \neq 0\}$ at $P^0$ is
\[
 z_3-\nabla\Phi(x^{0}) \cdot z'=c,
 \qquad z'=(z_1,z_2).
\]
In particular, this plane does not pass through the origin.

After subtracting this tangent plane, the local defining
function of $S_{c}$ near $P^{0}$ is
\[
 \phi(u,v):= v^kU(u,v).
\]
On a sufficiently small neighborhood of $(0,0)$, $U$ has constant
sign. The analytic change of coordinates
\[
 (u,v)\longmapsto
 \bigl(u,\,v|U(u,v)|^{1/k}\bigr)
\]
therefore puts the function $\phi$ into the form
$\pm  v^k$. Its Newton height is $k$.
By the invariance of height under local coordinate changes
and invertible linear changes of ambient coordinates, we have
\[
 h(P^0,S_{c})=k
 \le M_\Phi,
 \qquad
 S_{c}
 =\{(x,c+ \Phi(x)):x\ne0\}.
\]
Thus the transversality and finite-type hypotheses of
\cite[Theorem~1.2]{IKM10} hold at every point of the fixed
annulus. That theorem then yields the local $L^p$ range $p>k$. Then (\ref{away0}) follows from
a partition of unity on the global surface and the fact that $M_{\Phi} \ge k$.

\textbf{(2)}  Next, we consider   small $\varepsilon$. Apply the invertible
linear transformation
\[
 (y_1,y_2,y_3)\longmapsto(y_1,y_2,y_3/c).
\]
The transformed surface is the graph of
\[
 1+\frac{\varepsilon}{|c|}\widetilde\Phi(x),
 \qquad
 \widetilde\Phi=\frac{c}{|c|}\Phi.
\]
At every point of $\operatorname{supp}\tilde{\psi}$, Lemma \ref{lem:contactcone} gives a direction $e$ and an integer
$2\le k\le M_\Phi$ for which
$\partial_e^k\widetilde\Phi\ne0$. We cover the compact set
$\operatorname{supp}\tilde{\psi}$ by finitely many sufficiently small
patches on which the corresponding derivative remains
nonzero, and choose a smooth partition of unity.

On each patch, we apply  Proposition  \ref{quoted:epsilon} with
principal function $\widetilde\Phi$ and perturbation equal
to zero. Taking the minimum of the finitely many
small-parameter thresholds and summing the patch estimates,
we obtain constants $\varepsilon_0>0$ and $C_p$ such that
(\ref{eq:cone-fixed-annulus}) holds.
\end{proof}


Now we are ready to prove the \(L^{p} \)-boundedness of \(\mathcal{M}_{c}\).
Suppose that \(\operatorname{supp}\psi \subset \{x\in \mathbb{R}^{2}: |x| < s_{0}\}\) for
some constant \(s_{0}>0\). We choose a radial function
\(\chi_0\in C_c^\infty(\mathbb R^2)\) that is non-increasing as a function of
\(|x|\) and satisfies
\[
0\le\chi_0\le1,\qquad
\chi_0(x)=1\ \text{if }|x|\le1,\qquad
\chi_0(x)=0\ \text{if }|x|\ge2.
\]
Set
\[
\chi(x)=\chi_0(x)-\chi_0(2x),
\qquad
s_j=2^{-j}s_0,
\qquad
\chi_j(x)=\chi(x/s_j)
\quad (j\geq 0).
\]
Then \(\chi_j\ge0\) and
\[
\operatorname{supp}\chi_j
\subset\{x:s_j/2\le |x|\le2s_j\}.
\]
Moreover,
\begin{equation}\label{eq:vertex-partition}
    \sum_{j=0}^{\infty}\chi_j(x)
    =1
    \quad\text{for all }x.
\end{equation}
Given $\epsilon_0$ sufficiently small, we choose the smallest positive integer \(J\) such that \( s_{J}  < \epsilon_{0}\). For \(j\ge J\), define
\[
\mathcal M_{c,j}f(y)
:=
\sup_{t>0}\left|
\int_{\mathbb R^2}
f\bigl(y-t(x,c+\Phi(x))\bigr)
\psi(x)\chi_j(x)\,dx
\right|.
\]
We also define
\[
\mathcal M^{0}_{c}f(y)
:=
\sup_{t>0}\left|
\int_{\mathbb R^2}
f\bigl(y-t(x,c+\Phi(x))\bigr)
\psi(x)  \biggl( \sum_{j=0}^{J-1}\chi_j(x)  \biggr)\,dx
\right|.
\]

By \eqref{eq:vertex-partition} and the triangle inequality,
\begin{equation}\label{eq:near-dyadic-bound}
	\mathcal M_cf(y)
	\le \sum_{j=J}^{\infty}\mathcal M_{c,j}f(y) + \mathcal M^{0}_{c}f(y).
\end{equation}
By Proposition \ref{prop:cone-fixed-annulus}-(1), the operator $\mathcal{M}_{c}^{0}$ is
$L^{p} $-bounded for every $p > M_{\Phi}$. Thus, it remains to prove that
\[
\sum_{j \ge J} \|\mathcal{M}_{c,j}\|_{L^{p}  \rightarrow L^{p} }< +\infty
\]
for every $p> M_{\Phi}$.

Fix \(j\ge J\). After scaling \(x\) by \(s_{j}\) and using
\(\Phi(s_j x)=s_j\Phi(x)\), we obtain
\[
\mathcal M_{c,j}f(y)
=
s_j^2\sup_{t>0}\left|
\int_{\mathbb R^2}
f\bigl(y-t(s_j x,c+s_j\Phi(x))\bigr)
\psi(s_j x)\chi(x)\,dx
\right|.
\]
We define the invertible  linear map \(A_j\)
\[
A_j(y_1,y_2,y_3)
:=
\left(\frac{y_1}{s_j},\frac{y_2}{s_j},
 y_3 \right).
\]
  Consequently,
\[
A_j(s_j x,c+s_j\Phi(x))
=
\bigl(x,c+s_j \Phi(x)\bigr).
\]

We  define the normalized maximal
operator
\[
\widetilde{\mathcal M}_jf(y)
:=
\sup_{t>0}\left|
\int_{\mathbb R^2}
f\bigl(y-t(x,c+s_j \Phi(x))\bigr)
\psi(s_j x)\chi(x)\,dx
\right|.
\]
With $f_j(z)=f(A_j^{-1}z)$,  we have the exact relation
\[
\mathcal M_{c,j}f(A_j^{-1}y)
=
s_j^2\widetilde{\mathcal M}_j f_j(y).
\]
Consequently, the $L^p $-boundness for $\widetilde{\mathcal{M}}_{j}$ yields the estimate
\begin{equation}\label{eq:translatedcone-normscaling}
\|\mathcal M_{c,j}f\|_{L^p }
\le s_j^2
\|\widetilde{\mathcal M}_j\|_{L^p \to L^p }
\|f\|_{L^p }.
\end{equation}
Since $s_{j}< \epsilon_{0}$, Proposition \ref{prop:cone-fixed-annulus}-(2) gives
\[
 \|\widetilde{\mathcal M}_j\|_{L^p \to L^p }
 \le
 C_p  s_{j}^{-1/p}.
\]
Combining this estimate with
\eqref{eq:translatedcone-normscaling},
we obtain
\begin{equation}\label{eq:shellbound}
 \|\mathcal M_{c,j}\|_{L^p \to L^p }
 \le
 C_p s_j^2s_j^{-1/p}.
\end{equation}
Finally,
\[
 \sum_{j=J}^{\infty}s_j^{2-1/p}
 =
 s_0^{2-1/p}
 \sum_{j=J}^{\infty}2^{-j(2-1/p)}
 <
 \infty
\]
since $p>M_\Phi\ge2$. This completes the proof.

\section*{Acknowledgments}
In the early stages of this research, Prof. Detlef M\"{u}ller shared with us his ideas on the classification of mixed homogeneous functions whose Hessian determinant is identically zero. Although we later found that the complete classification can be obtained from the results in \cite{GN16}, we would still like to express our sincere gratitude to Prof. Detlef M\"{u}ller.

\section*{AI usage}
During the course of this research, we used GPT 6.1 to help us locate the reference \cite{GN16}. Using the results in \cite{GN16}, we completed the classification of mixed homogeneous functions whose Hessian determinant is identically zero. During the manuscript preparation stage, we used GPT 6.1 to assist with English writing. All derivations and calculations in this paper were performed manually by the authors, and the authors take responsibility for the correctness of the paper.

\begin{flushleft}

\vspace{0.3cm}\textsc{Wenjuan Li\\
	School of Mathematics and Statistics, Northwestern Polytechnical University, 710129, Xi'an, People's Republic of China\\
	Email address:} liwj@nwpu.edu.cn

\vspace{0.3cm}\textsc{Huiju Wang (Corresponding author)\\
	School of Mathematics and Statistics, Henan University, 475000, Kaifeng, People's Republic of China\\
	Email address:} huijuwang@mail.nwpu.edu.cn
\end{flushleft}

\end{document}